\documentclass[11pt, a4paper]{amsart}
\usepackage[lmargin=30mm,rmargin=30mm,tmargin=25mm,bmargin=25mm]{geometry}
\renewcommand{\baselinestretch}{1.15}
\usepackage{amssymb,amsmath}
\usepackage{amsthm}
\usepackage{amsfonts}
\usepackage{latexsym}
\usepackage[pdftex]{graphicx}
\usepackage{enumerate}
\usepackage{color}
\usepackage{cite}
\usepackage{tikz}

\usepackage[latin1]{inputenc}

\newtheorem{thm}{Theorem}

\newtheorem{lem}[thm]{Lemma}

\theoremstyle{definition}
\newtheorem{nota}[thm]{Notation}

\newcommand{\fibo}{Fibonacci number}
\newcommand{\stab}{stability number}
\newcommand{\ts}{tree of stars}    
\newcommand{\tss}{trees of stars} 
\newcommand{\tos}{\ts}    
                    
\newcommand{\ct}{center}       
\newcommand{\ctree}{center-tree}
\newcommand{\mss}{maximum stable set}
\newcommand{\gr}{good rotation}
\newcommand{\C}{C}

\newcommand{\DEF}{\sl}    
\newcommand{\expo}{exposed center}
\renewcommand{\l}{\ell}

\newcommand{\R}{\mathbb R{}}

\newcommand{\F}{F}

\title{Trees with Given Stability Number and Minimum Number of Stable Sets}

\author{V\'eronique Bruy\`ere}
\address{
\newline V\'eronique Bruy\`ere, Hadrien M\'elot
\newline Institut d'Informatique
\newline Universit\'e de  Mons -- UMONS
\newline 20 Place du Parc
\newline B-7000-Mons, Belgium
\newline {\tt \{veronique.bruyere, hadrien.melot\}@umons.ac.be}}

\author{Gwena\"el Joret}
\address{
\newline Gwena\"el Joret
\newline D\'epartement d'Informatique
\newline Universit\'e Libre de Bruxelles -- ULB
\newline Boulevard du Triomphe CP 212
\newline B-1050-Bruxelles, Belgium
\newline {\tt gjoret@ulb.ac.be}}

\author{Hadrien M\'elot}

\thanks{G. Joret is a
Postdoctoral Researcher of the {\em Fonds
National de la Recherche Scientifique (F.R.S.--FNRS)}}

\date{}

\begin{document}

\maketitle
\begin{abstract}
We study the structure of trees minimizing their number of stable sets 
for given order $n$ and stability number $\alpha$. 
Our main result is that the edges of a non-trivial extremal tree can
be partitioned into $n-\alpha$ stars, each of size $\left \lceil \frac{n-1}{n-\alpha} \right \rceil$ 
or $\left \lfloor \frac{n-1}{n-\alpha} \right \rfloor$, so that
every vertex is included in at most two distinct stars, 
and the centers of these stars form a stable set of the tree.
\end{abstract}

\section{Introduction}

The number $F(G)$ of stable sets (or independent sets) of a graph $G$ 
was first considered by Prodinger and Tichy~\cite{Prodinger82}. They called it the {\em Fibonacci number}
of a graph, based on the observation that the number of stable sets in a path on $n$
vertices is exactly the $n+2$-th Fibonacci number.
The invariant $F(G)$ is also known as the {\em Merrifield-Simmons index} or {\em $\sigma$-index} 
of the graph $G$~\cite{MS89}.

There is a rich literature dealing with extremal questions regarding $F(G)$, and the
case where $G$ is a tree received much attention.
For instance, Heuberger and Wagner~\cite{Heuberger08} characterized
the $n$-vertex trees $T$ with maximum degree $\Delta$ maximizing $F(T)$, for given $n$ and $\Delta$.
Li, Zhao, and Gutman~\cite{Li05}, and independently Pedersen and Vestergaard~\cite{Pedersen07}
determined the $n$-vertex trees $T$ maximizing $F(T)$ when the diameter is fixed. 
Yu and Lv~\cite{YuLv07} considered similarly the case of trees with a fixed number of leaves.
The reader is referred to \cite{Knopfmacher07, LinLin95, Wang08} 
and the references therein for other results of this kind on trees.
Other classes of graphs have been considered as well; this includes
unicyclic graphs~\cite{LiZhu09, Pedersen05, Pedersen06, Wang08a},
bicyclic graphs~\cite{DengChenZhang08},
tricyclic graphs~\cite{ZhuLiTan10},
quasi-trees~\cite{LiLiJing09},
maximal outerplanar graphs~\cite{Alameddine98},
connected graphs~\cite{Bruyere09},
and bipartite $d$-regular graphs~\cite{Kahn01}. 

In this paper, we propose to revisit an old result of Prodinger and Tichy~\cite{Prodinger82}
which is well-known in this area: among all trees $T$ on $n\geq 2$ vertices,
the path $P_{n}$ and the star $K_{1,n-1}$ respectively minimizes and maximizes $F(T)$.
Observe that the stability number (the largest cardinality of a stable set) 
of $P_{n}$ is $\lceil n/2 \rceil$, while that of $K_{1,n-1}$ is $n-1$. 
Since the stability number $\alpha$ of every tree on $n$ vertices is between these two
extreme values, one may wonder which are the trees $T$ 
minimizing or maximizing $F(T)$ for fixed $\alpha$. 
In a previous contribution~\cite{Bruyere09}, a subset of the authors showed  
that the trees maximizing $F(T)$ have a rather simple structure: they are exactly the
trees that can be obtained from the Tur\'an graph\footnote{
The {\DEF Tur\'an graph} $T(n,\alpha)$ is the $n$-vertex graph consisting of
$\alpha$ vertex-disjoint cliques which are as balanced as possible.}
$T(n,\alpha)$ by adding $\alpha -1$ edges.

\begin{figure}
   \centering
   \includegraphics[scale=0.9]{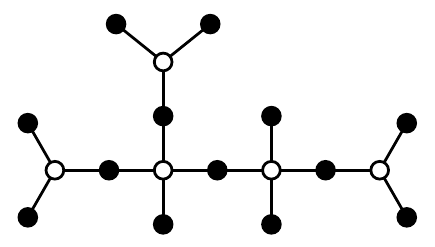} 
   \caption{\label{fig-intro}An extremal tree for $n= 18$ and $\alpha = 13$. 
   White vertices are the centers of the stars.}   
\end{figure}

Here, we focus on the trees minimizing $F(T)$ for fixed $\alpha$.
As will become apparent in the next sections, the structure of the 
corresponding extremal trees is less straightforward than when maximizing $F(T)$.
Our main result is that the edges of a (non-trivial) extremal tree can uniquely
be partitioned into $n-\alpha$ stars, each of size $\left \lceil \frac{n-1}{n-\alpha} \right \rceil$ 
or $\left \lfloor \frac{n-1}{n-\alpha} \right \rfloor$, so that
every vertex is included in at most two distinct stars, 
and the centers of these stars form a stable set of the tree.
(See Figure~\ref{fig-intro} for an illustration.)

It follows that an extremal tree can be completely described by specifying which 
stars have a vertex in common. This is captured by the following labeled tree,
which we call {\em \ctree}:
the tree has $n-\alpha$ vertices (representing the stars); 
there is an edge between two vertices if the corresponding stars have a vertex in common, 
and each vertex is labeled by the size of the corresponding star.
However, characterizing the {\ctree s} of extremal trees appears to be challenging.
We obtained such a characterization in a very restricted case, namely when
$(n - 1) \bmod (n-\alpha) \in \{0, n-\alpha - 2, n-\alpha - 1\}$, but
the general case is far from understood and is left as an open problem.

The paper is organized as follows. In Section~\ref{sec_prelim}, we provide the necessary 
definitions and a few basic lemmas. Sections~\ref{sec-task1} and~\ref{sec-task2},
which constitute the bulk of this paper, 
are devoted to the proof that extremal trees have a partition into stars as described above
(cf.\ Theorem~\ref{th-task2}).
Finally, in Section~\ref{sec-task3}, we conclude with some remarks on
the problem of characterizing the {\ctree s} of extremal trees, and provide a solution for the 
particular case mentioned above.

\section{Definitions and Preliminaries} \label{sec_prelim}

This section is devoted to basic definitions and notations used throughout the text.
Also, we introduce the notion of {\ts}, and prove some useful facts about these trees.

All graphs are assumed to be finite, simple, and undirected.
We generally follow the terminology and notations of Diestel~\cite{Diestel05}.
The neighborhood of a vertex $v$ of a graph $G$ is
denoted $N_{G}(v)$; also, we write $N_{G}[v]$ for the set $N_{G}(v) \cup \{v\}$.
The degree of $v$ is denoted by $\deg_{G}(v)$.
(We often drop the subscript $G$ when the graph is clear from the context.) 
We simply write $|G|$ for the order $|V(G)|$ of $G$.
A subset $S\subseteq V(G)$ of vertices of a graph $G$ is a
{\DEF stable set} if the vertices in $S$ are pairwise non adjacent in $G$.
The maximum cardinality of a stable set in $G$ is the {\DEF stability
number} of $G$, and is denoted $\alpha(G)$.

Recall that $\F(G)$ denotes the number 
of stable sets of $G$ (including the empty set).
In this paper, we call this invariant the {\DEF \fibo} of $G$, as in~\cite{Prodinger82}. 
The following easy properties  (see \cite{Gutman86, Li05, Prodinger82}) 
of the {\fibo} are used throughout the paper:

\begin{lem} \label{lem:basic} 
Let $G$ be a graph. If $G$ is not empty, then,
\begin{itemize}
\item $\F(G) = \F(G -v) + \F(G - N[v])$ for every vertex $v$ of $G$, and
\item $\F(G) = \prod_{i=1}^k \F(G_i)$
if $G$ is the disjoint union of $k$ graphs $G_1, G_{2}, \dots,  G_{k}$.
\end{itemize}

In particular,
\begin{itemize}
\item $\F(G) > \F(G -v)$ for every vertex $v$ of $G$, and
\item $\F(G) \le 2 \F(G-v)$,
with strict inequality if $\deg_{G}(v) \geq 1$.
\end{itemize}
\end{lem}

A tree $T$ is {\DEF extremal} if $T$ has minimum {\fibo} among 
all trees on $|T|$ vertices with {\stab} $\alpha(T)$. 
As mentioned in the introduction, 
the structure of these extremal trees is the main topic of this paper;
in particular, we will see that non-trivial extremal trees are {\tss} that are balanced 
(see below for the definitions).

In this paper, a {\DEF leaf} of a tree is defined as a vertex of degree  at most $1$. Thus in particular if the tree has a unique vertex, then this vertex is considered to be a leaf. (Usually leaves are required to have degree exactly $1$, but this definition will be more convenient for our purposes.).
A {\DEF \ts} is defined inductively as follows:
\begin{itemize}
\item a single vertex is a {\ts};
\item if $T_1, T_2, \ldots, T_k$ ($k \ge 2$) are disjoint {\tss} and $v_i$ is a leaf of $T_i$ for $i=1, 2, \ldots, k$,  
then the tree obtained from $T_{1} \cup \cdots \cup T_{k}$ by adding a new vertex $w$ adjacent to every vertex $v_{i}$ is also a {\ts}. 
\end{itemize}

\begin{figure}[!htb]
\begin{center}
\includegraphics[scale=0.9]{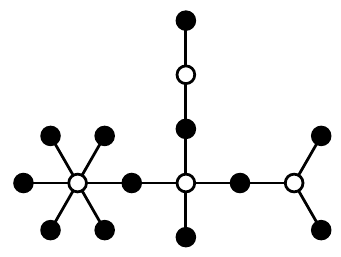}
\caption{\label{fig-exa-tos} 
A {\tos}. The white vertices are the centers of the tree.
}
\end{center}
\end{figure}

From the inductive definition given above, it is easy to check that, in a {\ts} $T$,
the distances between a fixed vertex $v\in V(T)$ and all the leaves of $T$ have the same
parity. The vertex $v$ is a {\DEF center} of $T$ is these distances are odd.
An example of a {\ts} is given in Figure~\ref{fig-exa-tos}. 

We note that {\tss} can equivalently be defined as follows. First observe that the vertex set of a tree $T$ can uniquely be partitioned into two stable sets $A$ and $B$. Then $T$ is a {\ts} if and only if one of these two sets, say $A$, contains all leaves of $T$ and no vertex with degree at least $3$ (thus all vertices in $A$ have degree at most $2$). In that case, the centers of $T$ are exactly the vertices in $B$. (We leave it to the reader to check that this definition of {\tss} is indeed equivalent to the original one.)

A {\ts} is {\DEF balanced} if the degrees of any two of its centers differ by at most~$1$ 
(for instance, the {\ts} in Figure~\ref{fig-intro} 
is balanced while the one in Figure~\ref{fig-exa-tos} is not).

Let $T$ be a {\ts}. The set of {\ct s} of $T$ is denoted by $\C(T)$.  
Observe that a vertex of $T$ with degree at least 3 is always a {\ct}; thus, the set 
$V(T) - \C(T)$ includes only vertices with degree at most 2. Also,
notice that $(\C(T), V(T) - \C(T))$ is a partition of $V(T)$ into two stable sets. 
Finally, note that a path $P$ is a {\ts} if and only if $P$ is even\footnote{
$P$ is odd (even) if $P$ has odd (respectively even) length, where the length
is defined as the number of edges.}.

\begin{lem}
\label{lem-mss}
Let $T$ be a {\ts}. Then $V(T) - \C(T)$ is the unique maximum stable set of $T$.
\end{lem}
\begin{proof}
We use induction on $|\C(T)|$. If $|\C(T)| = 0$, then $T$ is a single vertex 
and the claim trivially holds. Suppose $|\C(T)| \ge 1$, and
let $S$ be any {\mss} of $T$. Also, let $v$ be a leaf of $T$ 
and $w$ be its unique neighbor. 
Let $T_{1}, \dots, T_{k}$ be the components of $T - \{v, w\}$. Notice that each 
$T_{i}$  is a {\ts} and $\C(T) = \{w\} \cup \C(T_{1}) \cup \cdots \cup \C(T_{k})$.

First, suppose that $S\cap V(T_{i}) \neq  V(T_{i}) - \C(T_{i})$
for some $i\in \{1,\dots,k\}$. Then, by induction, $|S\cap V(T_{i})| < |V(T_{i}) - \C(T_{i})|$,
and the set 
$$
S' := \big(S - \left(S\cap V(T_{i})\right)\big) \cup \big(V(T_{i}) - \C(T_{i})\big)
$$
has cardinality larger than $S$. Since $S'$ cannot be a stable set of $T$, it follows $w \in S'$,
and thus $v \notin S'$.
But then $(S' - \{w\}) \cup \{v\}$ is a stable set of $T$, with cardinality equal to 
that of $S'$, a contradiction.

Therefore, $S\cap V(T_{i})=V(T_{i}) - \C(T_{i})$ for
every $i\in \{1,\dots,k\}$, which implies $w \notin S$, and hence $v\in S$. It follows that $S=V(T) - \C(T)$.
\end{proof}

\begin{figure}[!htb]
\begin{center}
\includegraphics[scale=0.9]{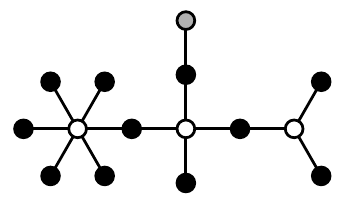}
\caption{\label{fig-exa-ats} 
A tree which is almost a {\tos}. Its {\expo} is drawn in grey.}
\end{center}
\end{figure}

Suppose that $T$ is a tree not isomorphic to a path and that $T$ has 
a leaf $w$ such that the tree $T'$ obtained by adding a new vertex $v$ adjacent to $w$ is a {\tos}.
Then $T$ is said to be {\DEF almost a {\tos}}, and
a vertex of $T$ is considered to be a center of $T$ if it is a center of $T'$.
Note that the requirement that $T$ is not a path ensures that the leaf $w$ is uniquely determined.
Hence, the set of centers of $T$ is well-defined.
The vertex $w$ is the unique center of $T$ which is also a leaf, and is said to be the {\DEF \expo} of $T$.
See Figure~\ref{fig-exa-ats} for an illustration.

Let us observe that, if $v$ is a leaf of a {\tos} $T$ with $|T| \geq 3$, then $T-v$ is exactly one of the following:
\begin{itemize}
\item a {\tos}; 
\item almost a {\tos};
\item an odd path.
\end{itemize}
This fact will be repeatedly used in our proofs.

\section{Extremal Trees are Trees of Stars}
\label{sec-task1}

The following theorem is a first step towards understanding the structure of extremal trees.

\begin{thm}
\label{th-task1}
Every extremal tree is either a {\ts} or an odd path.
\end{thm}

The proof of Theorem~\ref{th-task1} is based on the concept of rotating an edge in a tree $T$, which consists in first removing some edge $uv$ of $T$, and then adding an edge $uv'$, where $v'$ is some vertex of the component of $T-uv$ containing $v$.  
We will show that if $T$ is neither a {\ts} nor an odd path, then
there exists an edge of $T$ that can be rotated in such a way that the resulting tree $T'$
satisfies both $\alpha(T')=\alpha(T)$ and $F(T') < F(T)$.
To this aim, we introduce a few technical lemmas.

\subsection{Trees of Stars and the Golden Ratio}

The golden ratio $\phi := \frac{1 + \sqrt 5}{2}\simeq 1.618$ 
appears naturally when considering {\tss}, 
as illustrated by the following two lemmas.

\begin{lem}
\label{lem-ToS-leaf}
Let $T$ be a {\ts} and let $v$ be a leaf of $T$. Then 
$F(T) > \phi \cdot F(T-v)$. Moreover, this remains true if $T$ is almost a {\ts}, provided $v$ is not the {\expo}.
\end{lem}
\begin{proof}
The proof is by induction on $|T|$. The claim is true when $|T|=1$, since then
$F(T)=2>\phi \cdot F(T-v)=\phi$. 

For the inductive step, assume $|T| > 1$. Let $w$ be the center adjacent to $v$ 
and  $T_{1}, \dots, T_{k}$ be the components of $T-\{v,w\}$ (thus $k\geq 1$). 
Let also $v_{i}$ be the leaf of $T_{i}$ that is adjacent to $w$ in $T$. Letting
$$
\gamma := \prod_{i=1}^{k} \frac{F(T_{i})}{F(T_{i} - v_{i})},
$$
we deduce 
$$
\frac{F(T)}{F(T-v)} = \frac{2\prod_{i=1}^{k} F(T_{i}) 
+ \prod_{i=1}^{k} F(T_{i}-v_{i})}{\prod_{i=1}^{k} F(T_{i}) 
+ \prod_{i=1}^{k} F(T_{i}-v_{i})} = 2 - \frac{1}{1 + \gamma}.
$$

If $T$ is a {\ts}, then $T_{i}$ is also a {\ts} and $F(T_i) > \phi \cdot F(T_i-v_i)$ by the induction hypothesis,
for every $i\in\{1, \dots, k\}$. Hence, $\gamma >\phi^{k} \geq \phi$.

If, on the other hand,
$T$ is almost a {\ts}, then all trees $T_i$ are {\tss} except for exactly one, say $T_{1}$.
Thus, by induction, $F(T_i) > \phi \cdot F(T_i-v_i)$ for every $i\in\{2, \dots, k\}$.
Moreover, $T_{1}$ is either almost a {\ts} or an odd path. In the first case, $v_{1}$ is not the {\expo} of $T_{1}$, and
thus $F(T_1) > \phi \cdot F(T_1-v_1)$ by induction, which again implies $\gamma >\phi^{k} \geq \phi$.
In the second case, we have $k\geq 2$ since $T$ is not a path, and
$$
\gamma \geq \prod_{i=2}^{k} \frac{F(T_{i})}{F(T_{i} - v_{i})} > \phi^{k-1} \geq \phi, 
$$
by induction. 

Therefore, $\gamma > \phi$ holds in all possible cases, and it follows
$$
\frac{F(T)}{F(T-v)} =   2 - \frac{1}{1 + \gamma}> 
2 - \frac{1}{1 + \phi} = \phi.
$$
(The last equality is derived using the fact $\phi^{2} = \phi + 1$.)
\end{proof}

A similar but opposite property holds for centers:

\begin{lem}
\label{lem-ToS-center}
Let $T$ be a {\ts} and let $w$ be a center of $T$. Then 
$F(T) < \phi \cdot F(T-w)$. This moreover remains true if $T$ is almost a {\tos}.
\end{lem}
Note that, contrary to Lemma~\ref{lem-ToS-leaf}, here the vertex $w$ is allowed to be the {\expo} of $T$ if $T$ is almost a {\tos}.
\begin{proof}[Proof of Lemma~\ref{lem-ToS-center}]
Let  $T_{1}, \dots, T_{k}$ be the components of $T-w$. (Thus $k\geq 1$ and $k=1$ if and only if $w$ is the {\expo} of $T$.) 
We denote by $v_i$ the leaf of $T_i$ that is adjacent to $w$. Let
$$
\gamma := \prod_{i=1}^{k} \frac{F(T_{i})}{F(T_{i} - v_{i})}.
$$
We have
$$
\frac{F(T)}{F(T-w)} = \frac{\prod_{i=1}^{k} F(T_{i}) 
+ \prod_{i=1}^{k} F(T_{i}-v_{i})}{\prod_{i=1}^{k} F(T_{i})} = 1 + \frac{1}{\gamma}.
$$

If $T$ is a {\ts}, then $T_i$ is a {\ts} for every $i\in \{1, \dots, k\}$. The same is true if $T$ is almost a {\ts}
and $w$ is the {\expo} of $T$. Thus, in these two cases, 
we have $F(T_i) > \phi \cdot F(T_i-v_i)$ for every $i\in \{1, \dots, k\}$ by Lemma~\ref{lem-ToS-leaf}, 
which implies $\gamma > \phi^{k} \geq \phi$. 

Now, if $T$ is almost a {\ts} and $w$ is {\em not} the {\expo} of $T$, then $k\geq 2$ and all trees $T_i$ are {\tss} except for exactly one, 
say $T_{1}$. Thus, $F(T_i) > \phi \cdot F(T_i-v_i)$ for every $i\in\{2, \dots, k\}$ by Lemma~\ref{lem-ToS-leaf}.
Since $F(T_{1}) \geq F(T_1-v_1)$, this implies 
$$
\gamma \geq \prod_{i=2}^{k} \frac{F(T_{i})}{F(T_{i} - v_{i})} > \phi^{k-1} \geq \phi.
$$
It follows that $\gamma > \phi$ always holds, implying
$$
\frac{F(T)}{F(T-w)} =  1 + \frac{1}{\gamma} < 
1 + \frac{1}{\phi} = \phi.
$$
\end{proof}

\subsection{Edge Rotations}
\label{sec-rotations}

Let $T$ be a tree, $xy$ be one of its edges, and
$x'$ be a vertex distinct from $x$ in the component of $T -xy$ containing $x$.
Then the pair $\rho:=(yx, yx')$ is called a
{\DEF rotation}, and we let $\rho(T):= (T - yx) +yx'$ denote the 
tree resulting from the rotation.
The rotation $\rho$ is {\DEF good} if $\alpha(\rho(T)) = \alpha(T)$ and $F(\rho(T)) < F(T)$. 
(Thus, in order to prove that a tree is {\em not} extremal, it is enough to show
that it admits a good rotation.) 

When the rotation $\rho=(yx, yx')$ is clear from the context, we use the following notations.

\begin{nota}
\label{nota-rot}
Denote by $z$ ($z'$) the neighbor of $x$ (respectively $x'$) that
is included in the unique $xx'$-path in $T$. 
(Observe that $z=x'$ and $z'=x$ if $xx' \in E(T)$; similarly, $z=z'$
if $x$ and $x'$ are at distance $2$.) 
Let $x_{1}, \dots, x_{\l}$ be the neighbors of $x$ distinct from $z$ and $y$, 
and $x'_{1}, \dots, x'_{\ell'}$ be the neighbors of $x'$ distinct from $z'$.
These notations are illustrated in Figure~\ref{fig-rotating-edge}. 
Notice that $x$ and $x'$ have degree $\l+2$ and $\ell' +1$, respectively.  

\begin{figure}[!ht]
\begin{center}
\includegraphics{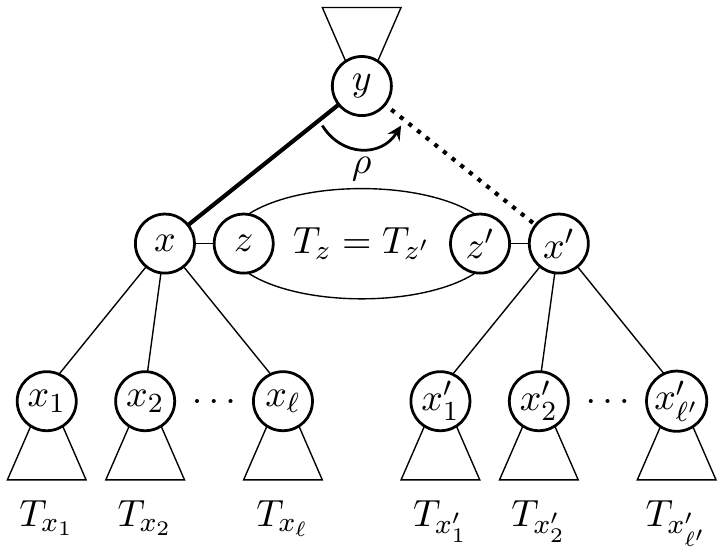}
\caption{\label{fig-rotating-edge} 
A rotation $\rho=(yx, yx')$. 
}
\end{center}
\end{figure}

Furthermore, we refer to the components of $T - \{ x,x'\}$ as follows:
If $v \in V(T) - \{x, x'\}$, then $T_{v}$ denotes the component of $T - \{ x,x'\}$
that includes $v$. Also, it will be convenient to define $T_{v}$ as the empty tree if $v\in\{x, x'\}$.
(In particular, $T_{z}=T_{z'}$, and the latter tree is empty if and only $xx' \in E(T)$.)
Finally, we define ten numbers corresponding to these components:
$$
\begin{array}{l l}
X  :=  \prod_{i=1}^{\l} F(T_{x_i}), &
\overline X  :=  \prod_{i=1}^{\l} F(T_{x_i} - x_{i}),\\[7pt]
X'  :=  \prod_{i=1}^{\ell'} F(T_{x'_i}), & 
\overline X'  :=  \prod_{i=1}^{\ell'} F(T_{x'_i} - x'_{i}),\\[7pt]
Y  :=  F\left(T_{y} \right), &
\overline Y := F \left(T_{y} - y \right),\\[7pt]
Z :=  F \left(T_{z} \right), &
\overline Z_{z} :=  F \left(T_{z} - z \right),\\[7pt]
\overline Z_{z'} :=  F \left(T_{z} - z' \right) , \textrm{ and} & \overline Z_{\{z,z'\}} :=  F \left(T_{z} - {\{z,z'\}} \right).
\end{array}
$$
Since the {\fibo} of an empty graph is $1$, by convention we let $X = \overline X := 1$ if $\l = 0$. Similarly, 
$X' = \overline X' := 1$ if $\ell' = 0$, and $Z = \overline Z_{z} =  \overline Z_{z'} = \overline Z_{\{z,z'\}} := 1$
if $x$ is adjacent to $x'$.
\end{nota}

\begin{lem}
\label{lem-rotating-edge}
With Notation~\ref{nota-rot}, we have
$F(\rho(T)) < F(T)$ if and only if $X \overline X' \overline Z_{z'} > \overline X X' \overline Z_{z}$. 
\end{lem}
\begin{proof}
Let us compute $F(T)$ by applying twice Lemma~\ref{lem:basic} (first item) with vertices $x$ and $x'$.  
We obtain
$$
F(T) = X  Y \left( X' Z + \overline X' \ \overline Z_{z'} \right) + {\overline X} \ {\overline Y}  \left( X' \overline Z_{z} + \overline X' \ \overline Z_{\{z,z'\}} \right).
$$
Similarly,
$$
F(\rho(T)) =  X'  Y \left( X Z + \overline X \ \overline Z_{z} \right) + {\overline X}' \ {\overline Y}  \left( X \overline Z_{z'} + \overline X \ \overline Z_{\{z,z'\}} \right). 
$$
Hence, we get
\begin{align*}
F(T) - F(\rho(T)) &= X \overline X' \overline Z_{z'} \left( Y - \overline Y \right) + \overline X X' \overline Z_{z} \left(\overline Y - Y \right),\\
&= \left(X \overline X' \overline Z_{z'} - \overline X X' \overline Z_{z} \right) \cdot \left(Y - \overline Y\right).
\end{align*}
Since $Y - \overline Y > 0$, we deduce that
$F(\rho(T)) < F(T)$ if and only if $X \overline X' \overline Z_{z'} > \overline X X' \overline Z_{z}$. 
\end{proof}

For the remainder of this section, let $T$ be a tree with a vertex $v$ of degree $k + 1 \geq 3$ such that the components of $T - v$ can be denoted as $T', T_{1}, T_{2}, \dots, T_{k}$ in such a way that, 
letting $v_{i}$ ($i\in \{1, \dots, k\}$) be the neighbor of $v$ in $V(T_{i})$ 
and letting $T^{+}_{i}$ be the tree obtained from $T_i$ by adding the vertex $v$ and the edge $vv_i$, the following holds:
\begin{itemize}
\item $T^{+}_{1}$ is a {\tos}, and
\item for each $i\in \{2, \dots, k\}$, at least one of the following two conditions holds:
\begin{enumerate}[(C1)]
\item \label{type1} $T_{i}$ is a {\tos} and $v_{i}$ is a leaf of $T_{i}$,
\item \label{type2} $T^{+}_{i}$ is a {\tos}.
\end{enumerate}
\end{itemize}
Also, let $v'$ be the neighbor of $v$ in $V(T')$, and let
$w$ be a leaf of $T_{1}$ distinct from the vertex $v_{1}$ (such a leaf exists
since $|T^{+}_{1}| \geq 3$).
Finally, let us emphasize that no assumption is made on the component $T'$, that is, $T'$
is an arbitrary non-empty tree. 

The following lemma is a crucial tool in our proof of Theorem~\ref{th-task1}:

\begin{lem} \label{lem-gr}
Let $T$ be as above. Then there exists a {\gr} in $T$.
\end{lem}

In order to prove Lemma~\ref{lem-gr}, we distinguish three cases:
\begin{itemize}
\item $T_{1}$ is a path and (C\ref{type1}) holds for some $i\in\{2,\dots, k\}$;
\item $T_{1}$ is not a path and (C\ref{type1}) holds for some $i\in\{2,\dots, k\}$, and
\item (C\ref{type2}) holds for every $i\in\{2,\dots, k\}$.
\end{itemize}
(Observe that $T$ falls in at least one of these cases.)
We consider each of these cases separately;
Lemma~\ref{lem-gr} is obtained by combining
Lemmas~\ref{lem-case1-gr}, \ref{lem-case2-gr} and \ref{lem-case3-gr}, which
respectively address the first, second, and third case.
The latter lemmas rely in turn on Lemmas~\ref{lem-case1_2-mss} and~\ref{lem-case3-mss}
below, showing that some specific rotations do not change the stability number of $T$.

\begin{figure}
\begin{center}
\includegraphics{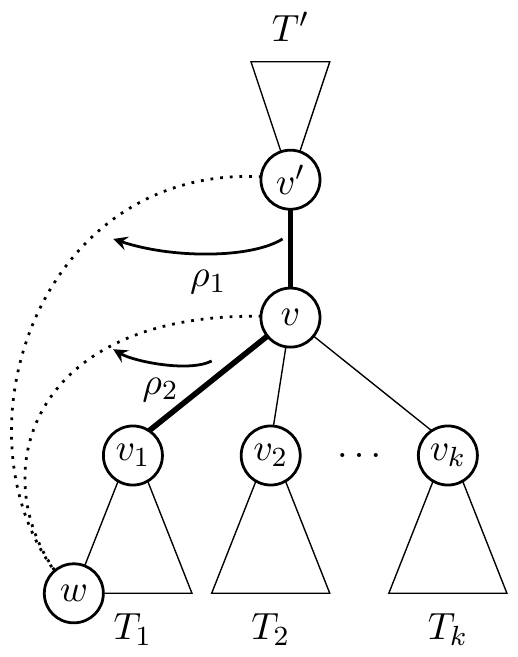}
\caption{\label{fig-two-rot} 
Rotations $\rho_1 = (v'v, v'w)$ and $\rho_2 = (vv_1, vw)$.
}
\end{center}
\end{figure}

Let $\rho_1$ and $\rho_2$ denote the two rotations $\rho_1 := (v'v, v'w)$ and $\rho_2 := (vv_1, vw)$
(see Figure~\ref{fig-two-rot} for an illustration).

\begin{lem} \label{lem-case1_2-mss}
Suppose that (C\ref{type1}) holds for some $i\in\{2,\dots, k\}$.
Then $\alpha(T)=\alpha(\rho_2(T))$, and $\alpha(T)=\alpha(\rho_1(T))$ if $T_{1}$ is a path.
\end{lem}
\begin{proof}
Assume without loss of generality that (C\ref{type1}) holds for $i=2$.
We first show:
\begin{equation}
\label{eq:avoid-v}
\textrm{Each of $T$, $\rho_1(T)$, and $\rho_2(T)$ contains a {\mss} not including $v$.}
\end{equation}
First, consider  the tree $T$. Let $S$ be a {\mss} of $T$. Suppose that $v \in S$ (otherwise, we are trivially done). Then $v_{2} \notin S$.
Since $T_2$ is a {\tos}, $T_{2}$ has a unique {\mss}, namely $S_2:=V(T_{2}) - \C(T_{2})$ 
(see Lemma~\ref{lem-mss}). 
In particular, $v_{2} \in S_{2}$, implying $| S \cap V(T_2) | < | S_2 |$. Let 
$$
S' := \big(S - (V(T_2)\cup \{v\})\big) \cup S_2.
$$
Observe that $S'$ is also a stable set of $T$, and $|S'| \geq |S|$. It follows that $S'$ is a {\mss} of $T$ with $v\notin S'$, as claimed.
Finally, note that the above argument still holds if we replace the tree $T$ by $\rho_1(T)$ or $\rho_2(T)$, 
which completes the proof of \eqref{eq:avoid-v}.

By~\eqref{eq:avoid-v}, there is a {\mss} $S$ of $T$ such that $v\notin S$.
The set $S$ is also a stable set of $\rho_2(T)$, implying $\alpha(\rho_2(T)) \ge \alpha(T)$.
Conversely, there is a {\mss} $S'$ of $\rho_2(T)$ with $v\notin S'$.
Since $S'$ is a stable set of $T$, this shows $\alpha(\rho_2(T)) \le \alpha(T)$, and hence
$\alpha(T) = \alpha(\rho_2(T))$. 

Now, suppose that $T_{1}$ is a path. Observe that $T_{1}$ has odd length, since $T^+_1$ is a {\tos}.
By~\eqref{eq:avoid-v}, there is a {\mss} $S'$ of $\rho_1(T)$ such that $v\notin S'$.
This set is also a stable set of $T$, implying $\alpha(T) \ge \alpha(\rho_1(T))$.
Let $S$ be a {\mss} of $T$. If $v'\notin S$ or  $w \notin S$ then $S$ is also a stable set of $\rho_1(T)$, showing
$\alpha(T) \le \alpha(\rho_1(T))$, and hence $\alpha(T) = \alpha(\rho_1(T))$. Assume thus $v',w \in S$. 
Then $v \notin S$. Since $T_{1}$ is a path with an even number of vertices, $V(T_{1})$ can be uniquely partitioned into two maximum stable sets $S_{1}$ and $S_{2}$ of $T_{1}$. Assuming without loss
of generality $w\in S_{1}$, we then have $S \cap V(T_1) = S_{1}$, because otherwise we could modify $S$
and find a larger stable set of $T$. It follows that
$$
\tilde S := (S - S_{1}) \cup S_{2}
$$
is another stable set of $T$ with $|\tilde S|=|S|$ and $w\notin \tilde S$, and we deduce
that $\alpha(T) = \alpha(\rho_1(T))$, as previously. 
\end{proof}

\begin{lem} \label{lem-case3-mss}
Suppose that (C\ref{type2}) holds for every $i\in\{2,\dots, k\}$. Then $\alpha(T) = \alpha(\rho_1(T))$. 
\end{lem}
\begin{proof}
First, we show that both $T$ and $\rho_1(T)$ contain a {\mss} that does not include $v'$. 
Let $S$ be a {\mss} of $T$ (respectively, $\rho_1(T)$), and suppose that $v' \in S$. Then, $v \notin S$ (respectively, $w \notin S$).
The tree $T^+_i$ ($i\in\{1,\dots, k\}$) is a {\tos}, and hence 
contains a unique {\mss} $S^+_i$ by Lemma~\ref{lem-mss}. (Note that $v, w \in S^+_1$, and
$v \in S^+_i$ for every $i\in\{2,\dots, k\}$.)
Since $v \notin S$ or $w \notin S$, it follows that $| S \cap V(T^+_1) | <  | S^+_1 |$. Also,
$| S \cap V(T^+_i) |  \le  | S^+_i |$ for every $i\in\{2,\dots, k\}$.
Let
$$
S' := \big(S \cap (V(T') - \{v'\})\big) \cup S^+_1 \cup S^+_2 \cdots \cup S^+_k.
$$
The set $S'$ is a stable set of $T$ (respectively, $\rho_1(T)$), and $|S'| \geq |S|$ by the previous observations, implying that 
$S'$ is a {\mss} with $v' \notin S'$.

Now, observe that every stable set of $T$ (respectively $\rho_1(T)$) that does not include $v'$ is also a stable set of $\rho_1(T)$
(respectively $T$). Since there exists a {\mss} not including $v'$ in each of these two trees, it follows $\alpha(T) = \alpha(\rho_1(T))$.
\end{proof}

\begin{lem} \label{lem-case1-gr}
Suppose that $T_{1}$ is a path and that (C\ref{type1}) holds for some $i\in\{2,\dots, k\}$.
Then the rotation $\rho_1$ is good for $T$. 
\end{lem}
\begin{proof}
We have that $\alpha(T) = \alpha(\rho_1(T))$ by Lemma~\ref{lem-case1_2-mss}, thus
it remains to prove that $F(\rho_1(T)) < F(T)$. 
Let $x:=v, x':=w, y:=v'$, and consider Notation~\ref{nota-rot} 
with respect to the rotation $\rho_1=(v'v, v'w)=(yx,yx')$. 
We have $X'=\overline X' = 1$ and $\overline Z_z = \overline Z_{z'}$ since $T_{1}$ is a path.
Also, $\overline X < X$, because $\ell=k-1 > 0$. Therefore, 
$X \overline X' \overline Z_{z'} > \overline X X' \overline Z_{z}$, and hence $F(\rho_1(T)) < F(T)$
by Lemma~\ref{lem-rotating-edge}.
\end{proof}

\begin{lem} \label{lem-case2-gr}
Suppose that $T_{1}$ is not a path and that (C\ref{type1}) holds for some $i\in\{2,\dots, k\}$.
Then the rotation $\rho_2$ is good for $T$. 
\end{lem}
\begin{proof}
Since (C\ref{type1}) holds for some $i\in\{2,\dots, k\}$, 
Lemma~\ref{lem-case1_2-mss} implies that $\alpha(T) = \alpha(\rho_2(T))$.  
We prove that $F(\rho_2(T)) < F(T)$, using Lemma~\ref{lem-rotating-edge} again. 
Let $x:=v_{1}, x':=w, y:=v$, and consider the notations 
associated to the rotation $\rho_2=(vv_{1}, vw)=(yx,yx')$.
(Thus, in particular, $z$ and $z'$ are the neighbors of respectively $v_1$ and $w$ that are included in the unique $v_1w$-path.)
We have $X'=\overline X' = 1$ since $w$ is a leaf. Hence, by Lemma~\ref{lem-rotating-edge}, 
to prove $F(\rho_2(T)) < F(T)$ it suffices to show
that $\overline X \ \overline Z_z < X \overline Z_{z'}$.

If $v_{1}w \in E(T)$, then $T_z$ is empty and $\overline Z_z = \overline Z_{z'} = 1$. Moreover, $\l \ge1$ because $T_1$ is not a path,
implying $\overline X < X$ and $F(\rho_2(T)) < F(T)$.

Now, assume that $v_{1}w \notin E(T)$. 
Since $v_{1}$ is a center of the {\tos} $T^+_1$, its distance to the leaf $w$ is odd, and hence at least $3$. This implies $z \neq z'$.
Also, the component of $T^+_1-v_{1}$ that includes $w, z$, and $z'$ is also a {\tos}; hence,
$T_{z}$ is either a {\tos}, or almost a {\tos}, or an odd path.

In the first two cases,
$z$ and $z'$ are respectively a leaf and a center of $T_{z}$. Furthermore, $z'$ is the {\expo} of $T_{z}$
in the second case. It follows $Z > \phi \cdot \overline Z_z$ and $Z < \phi \cdot \overline Z_{z'}$
from respectively Lemmas~\ref{lem-ToS-leaf} and~\ref{lem-ToS-center}, implying $\overline Z_z < \overline Z_{z'}$.
This in turn implies $\overline X \ \overline Z_z < X \overline Z_{z'}$ since trivially $\overline X \le X$.

In the third case, we have that $\overline Z_z = \overline Z_{z'}$, because $T_{z} - z$ is isomorphic to 
$T_{z} - z'$. Since $T_{1}$ itself is not a path, we also have $\ell = \deg_{T}(v_{1}) - 2 > 0$, implying 
$\overline X < X$. It follows again $\overline X \ \overline Z_z < X \overline Z_{z'}$.
\end{proof}

\begin{lem} \label{lem-case3-gr}
Suppose that (C\ref{type2}) holds for every $i\in\{2,\dots, k\}$.
Then $T$ admits a {\gr}.
\end{lem}

\begin{proof}
The proof involves two different rotations. 
First, let $x:=v, x':=w, y:=v'$, and consider the rotation $\rho_{1}=(v'v, v'w)=(yx,yx')$ and the corresponding notations. 
(Thus, $z=v_{1}$ and $z'$ is the unique neighbor $w'$ of $w$.)
We have $X'=\overline X' = 1$ and $\overline X < X$. 

If $\overline Z_{z} \le \overline Z_{z'}$, then $F(\rho_1(T)) < F(T)$ by Lemma~\ref{lem-rotating-edge}. 
Also, $\alpha(T) = \alpha(\rho_1(T))$ by Lemma~\ref{lem-case3-mss}, implying that 
$\rho_{1}$ is a good rotation. 
Thus, we may assume that $\overline Z_{z} > \overline Z_{z'}$, that is,
\begin{equation}
\label{eq:T1}
F(T_1 - \{w, v_1\}) > F(T_1 - \{w, w'\}).
\end{equation}
(Observe that this implies $w' \neq v_{1}$.)

\begin{figure}
\begin{center}
\includegraphics[scale=0.9]{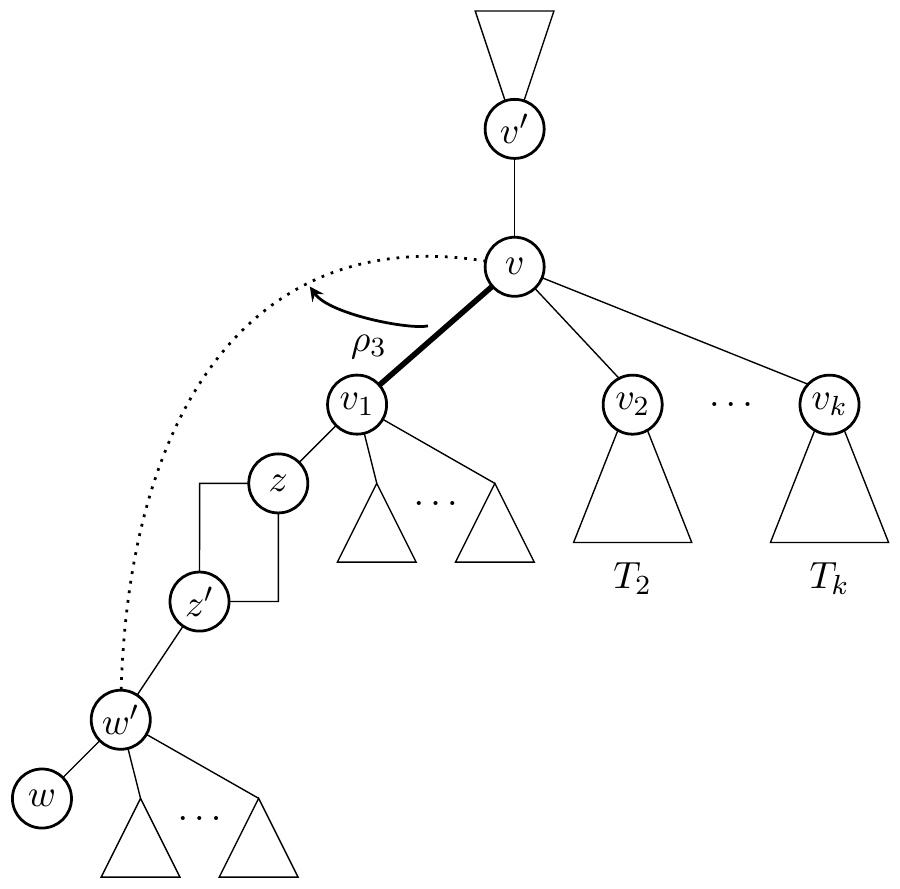}
\caption{\label{fig-second-rot} 
The rotation $\rho_3 = (vv_{1}, vw')$ in the proof of Lemma~\ref{lem-case3-gr} (note that possibly $z=z'$).
}
\end{center}
\end{figure}

Now, let $\rho_3$ be the rotation  $\rho_3:= (vv_1,vw')$ (see Figure~\ref{fig-second-rot} for an illustration).
We will show that $\rho_3$ is a good rotation.
Similarly as before, let $x:=v_{1}, x':=w', y:=v$, and consider the notations associated to $\rho_3 = (vv_1,vw')=(yx,yx')$.
We have
$$
F(T_1-\{w,v_1\}) = X(X'Z + \overline X' \overline Z_{z'})
$$
and
$$
F(T_1-\{w,w'\}) = X'(XZ + \overline X \ \overline Z_{z}).
$$
Thus,~\eqref{eq:T1} implies that $X \overline X' \overline Z_{z'} > \overline X X' \overline Z_{z}$, and it follows from
Lemma~\ref{lem-rotating-edge} that $F(\rho_3(T)) < F(T)$.
Therefore, it remains to show that $\alpha(T) = \alpha(\rho_3(T))$.

Let $S$ be a {\mss} of $T$. We may assume that $w' \notin S$. (Indeed, if not, consider the set $(S - \{w'\}) \cup \{w\}$ instead.)
The set $S$ is also a stable set of $\rho_3(T)$, showing $\alpha(T) \leq \alpha(\rho_3(T))$.

Now, let $S$ be a {\mss} of $\rho_3(T)$.
If $v\notin S$ or $v_{1}\notin S$, then $S$ is a stable set of $T$, showing $\alpha(T) \geq \alpha(\rho_3(T))$,
and hence $\alpha(T) = \alpha(\rho_3(T))$. Thus, suppose that $v, v_{1} \in S$.
The set $\tilde S:=(S - \{v_{1}\}) \cap V(T^{+}_{1})$ is a stable set of the {\tos} $T^{+}_{1}$. 
Moreover, $|\tilde S| < |S^{+}_{1}|$ since $z\notin \tilde S$, where $S^{+}_{1}:=V(T^{+}_{1}) - \C(T^{+}_{1})$ is
the unique {\mss} of $T^{+}_{1}$ (see Lemma~\ref{lem-mss}). This implies that the set
$$
(S - V(T^{+}_{1})) \cup S^{+}_{1}
$$
has cardinality at least that of $S$, and furthermore is a stable set of $T$. 
Again, it follows that $\alpha(\rho_3(T)) = \alpha(T)$.
\end{proof}

\subsection{Proof of Theorem~\ref{th-task1}}

Now, we may turn to the proof of Theorem~\ref{th-task1}.

\begin{proof}[Proof of Theorem~\ref{th-task1}]
Let $T$ be an extremal tree. Arguing by contradiction, assume that $T$ is neither a path nor a {\ts}. 
(Recall that even paths are {\tss}.)

Root the tree $T$ at an arbitrary leaf $r$.
Let $T_{v}$ ($v\in V(T)$) be the subtree of $T$ rooted at vertex $v$.
To each leaf $u$ of $T$ we associate a corresponding {\DEF witness}, defined
as the highest ancestor $v$ of $u$ such that $T_{v}$ is a {\tos} {\em and} $v$ is a 
leaf of $T_{v}$. Since $u$ itself satisfies these two conditions, the latter
vertex is well-defined. 

Let us look at a few properties of witnesses: First, clearly $r$ is not a witness (for otherwise
$T$ would be a {\tos}). Also, if $v$ and $w$ are two distinct witnesses, then
$v$ is neither an ancestor nor a descendant of $w$. (In other words, the set of
witnesses forms an antichain in the partial order implied the rooted tree.)

We may assume that $r$ has been chosen so that it satisfies:
\begin{equation}
\label{eq-change-root}
\textrm{
If the neighbor of $r$ has degree $2$ in $T$, then $T - r$ is not a {\ts}.
}
\end{equation}
Indeed, let $v$ be the neighbor of $r$, and
suppose $\deg_{T}(v)=2$ and that $T - r$ is a {\tos}. Since $T - r$ is not a path, 
there is a vertex $w$ with degree at least 3 in $T-r$, which is therefore a center of $T-r$.
Since $v$ is a leaf of $T-r$, it follows that the distance in $T$ between $r$ and $w$ is even.
Now, consider a component of $T-w$ that does not contain $r$, and select a leaf $r'$ of that component.
The tree $T-r'$ cannot be a {\ts}, because $w$ is at even distance from the leaf $r$ in $T-r'$.
Thus, we deduce that~\eqref{eq-change-root} holds if we root $T$ at $r'$ instead of $r$.

Let $u$ be a witness of maximum depth in $T$, and let $v_{1}$ be the parent of $u$.
It follows from~\eqref{eq-change-root} that $v_{1}\neq r$, hence $v_{1}$ has a parent $v$.
Let $T_{1}^{+}$ be the tree obtained from $T_{v_1}$ by adding the vertex $v$ and the edge $vv_1$. 
We show: 
\begin{equation}
\label{eq-T1}
\textrm{
$T_{1}^{+}$ is a {\tos}.
}
\end{equation}
This is trivially true if $\deg(v_{1})=2$, hence assume $\deg(v_{1})\geq 3$,
and let $u_{1}, \dots, u_{\ell}$ denote the children of $v_{1}$ distinct from $u$.
Let $j\in \{1, \dots, \ell\}$, and
consider a leaf $z$ of $T$ which is contained in $T_{u_{j}}$.
In $T$, the witness $w$ of $z$ cannot be higher than $u_{j}$ (since
otherwise $w$ would be an ancestor of the witness $u$). Also, $w$ cannot be
a descendant of $u_{j}$, because it would contradict the fact that $u$ has maximum
depth among all witnesses. Thus $w=u_{j}$, and hence $T_{u_{j}}$ is a {\tos}
and $u_{j}$ is a leaf of $T_{u_{j}}$. It follows that 
$T_{1}^{+}$ is also a {\tos}.

The vertex $v$ has at least two children (counting $v_{1}$), since
otherwise $T_{v} = T^{+}_{1}$, contradicting the fact that $u$ is
a witness. Observe that this implies $v \neq r$.
Let $v_{2}, \dots, v_{k}$ be the children of $v$ that are
distinct from $v_{1}$, and let $T_{i}^{+}$ ($i\in \{2, \dots, k\}$) 
be the tree obtained from $T_{v_i}$ by adding the vertex $v$ and the edge $vv_i$. 
 
Let $i\in \{2, \dots, k\}$, and consider any witness $u'$
contained in $V(T_{v_{i}})$ (observe that there must be at least one).
Since the depth of $u'$ is at most that of $u$, either $u'=v_{i}$ or $u'$ is a child of $v_{i}$.
In the first case, $T_{v_{i}}$ is a {\tos} and $v_{i}$ is a leaf of that tree, 
by definition of a witness.
In the second one, since $u$ and $u'$ have the same depth, the proof 
of~\eqref{eq-T1} directly shows that $T_{i}^{+}$ 
is a {\tos} (one just needs to replace $u$ by $u'$ and $v_{1}$ by $v_{i}$).

Let $T'$ be the component of $T - v$ containing the root $r$.
Also, let $T_{i} := T_{v_{i}}$ for $i \in \{1, \dots, k\}$.
Let us summarize the previous observations: 
$T_{1}^{+}$ is a {\tos}, and for every $i \in \{2, \dots, k\}$, either 
$T_i$ is a {\tos} and $v_{i}$ is a leaf of $T_{i}$, or $T_{i}^{+}$ is a {\tos}.
It follows that $T$ satisfies the requirements of Lemma~\ref{lem-gr}, and therefore
contains a good rotation, contradicting the fact that $T$ is extremal. 
\end{proof}

\section{Extremal Trees of Stars are Balanced}
\label{sec-task2}

We have seen that every extremal tree is a {\ts} or an odd path (cf.\ Theorem~\ref{th-task1}).
In this section, we refine this result by showing that,
if an extremal tree $T$ is a {\ts}, then $T$ must be balanced.
(Recall that a {\ts} is balanced if the degrees of every two of its centers differ by at most $1$.)

\begin{thm}
\label{th-task2}
Every extremal tree is either a balanced {\ts} or an odd path.
\end{thm}

We again resort to edge rotations to prove Theorem~\ref{th-task2}. To this aim,
we need to introduce a few additional lemmas.

For an integer $k\geq 2$, let $f_{k}: \R \to \R$ be the function defined as
$$
f_{k}(x) := x^{k} - x^{k-1} + 2x - 1.
$$

\begin{lem}
\label{lem-f-increasing}
The function $f_{k}(\cdot)$ is strictly increasing on the interval $[0,1]$.
\end{lem}
\begin{proof}
We prove that the first derivative $f'_{k}(\cdot)$ is strictly positive on $[0,1]$, 
by induction on $k$. Let $x\in [0,1]$. We have $f'_{k}(x) > 0$ when $k=2$, since $f'_{2}(x) = 2x + 1 > 0$. 
For the inductive step, let $k \geq 3$, 
and rewrite $f'_{k}(x)$ as follows:
\begin{align*}
f'_{k}(x) &= kx^{k-1} - (k-1)x^{k-2} + 2 \\
&= 2 + x\left( (k-1)x^{k-2} + x^{k-2} - (k-2)x^{k-3} - x^{k-3}  \right) \\
&= 2 + x\left( f'_{k-1}(x) - 2 \right) + x^{k-1} - x^{k-2}.
\end{align*}
Since $f'_{k-1}(x) - 2 > -2$ by the induction hypothesis, we deduce
\begin{align*}
f'_{k}(x) &> 2 - 2x + x^{k-1} - x^{k-2}\\
&=(1-x)(2-x^{k-2}) \\
&\ge 0,
\end{align*}
as claimed.
\end{proof}

Since $f_{k}(1/2)= -1/2^{k} < 0 < 1 = f_{k}(1)$ and $f_{k}(\cdot)$ is continuous, 
it follows from Lemma~\ref{lem-f-increasing}
that $f_{k}(\cdot)$ has a unique root in the open interval $(1/2, 1)$, which we denote by $R_{k}$. 

\begin{lem}
\label{lem-Rk}
$R_{k} < 2^{-k/(k+1)}$ for every $k\geq 2$.
\end{lem}
\begin{proof}
Let 
$$
g(k) := f_{k}\left(2^{-k/(k+1)}\right)
= 2^{- k^{2}/(k+1)} - 2^{ - k(k-1)/(k+1)} 
+ 2 \cdot 2^{- k/(k+1)} - 1.
$$
By Lemma~\ref{lem-f-increasing}, it is enough to show $g(k) > 0$ for every integer $k \geq 2$. 
A quick hand-on computation shows that this is true for $k=2$ and $3$. 
Hence, we may assume $k\geq 4$. 

We first rewrite $g(k)$ as follows:
\begin{align*}
g(k) &= 2^{- k^{2}/(k+1)} - 2^{ - k(k-1)/(k+1)} 
+ 2 \cdot 2^{- k/(k+1)} - 1 \\
&= 2^{- k^{2}/(k+1)} - 2^{ k/(k+1) -  k^{2}/(k+1)} 
+ 2^{1/(k+1)} - 1 \\
&=  2^{1/(k+1)} - 1 
- \frac{2^{k/(k+1)} - 1}{2^{k^{2}/(k+1)}}.
\end{align*}
Using  $k^{2}/(k+1)  > k - 1$ and $2^{k/(k+1)} < 2$, we then obtain 
\begin{align*}
g(k) &>   2^{1/(k+1)} - 1 - \frac{2^{k/(k+1)} - 1}{2^{k-1}} \\
&> 2^{1/(k+1)} - 1 - 2^{-(k-1)}.
\end{align*}

For every real $y>0$, we have $e^{y} > 1  + y$, and hence $2^{y} > 1 + y\ln 2$,
since $2^{y} = e^{y\ln2}$.  It follows
$$
g(k) > \frac{\ln 2}{k+1} - \frac{1}{2^{k-1} } 
= \frac{1}{k+1}\left(\ln 2 - \frac{k+1}{2^{k-1}} \right).
$$
Since the function $\ln 2 - (k+1)/(2^{k-1})$ is increasing in $k$ and positive for $k=4$,
we deduce that $g(k) > 0$, as claimed.
\end{proof}

\begin{lem}
\label{lem-ToS-leaf-gen}
Let $k\geq 2$ be an integer, let $T$ be a {\ts} such that all its centers have degree at least $k$, and
let $v$ be a leaf of $T$. Then $F(T - v) < R_{k} \cdot F(T)$.
\end{lem}
Observe that $R_{2}=1 / \phi$; hence, this lemma generalizes Lemma~\ref{lem-ToS-leaf} in the case of {\tss}.
\begin{proof}[Proof of Lemma~\ref{lem-ToS-leaf-gen}]
The proof is by induction on $|T|$. The claim is true when $|T|=1$, since then
$F(T-v) / F(T) = 1/2 < R_{k}$.

For the inductive step, assume $|T| > 1$, and
let $w$ be the unique neighbor of $v$. Each component of $T-\{v,w\}$ is a {\ts};
let us denote these trees by $T_{1}, \dots, T_{\ell}$ (thus $\ell\geq k-1$). Let also $v_{i}$
be the leaf of $T_{i}$ that is adjacent to $w$ in $T$. 

Letting
$$
\gamma := \prod_{i=1}^{\ell} \frac{F(T_{i} - v_{i})}{F(T_{i})},
$$
we deduce that
$$
\frac{F(T-v)}{F(T)} = \frac{\prod_{i=1}^{\ell} F(T_{i}) 
+ \prod_{i=1}^{\ell} F(T_{i}-v_{i})}{2\prod_{i=1}^{\ell} F(T_{i}) 
+ \prod_{i=1}^{\ell} F(T_{i}-v_{i})} =\frac{1 + \gamma}{2 + \gamma}.
$$
The induction hypothesis gives  
$$
0 < \gamma < (R_{k})^{\ell} \leq (R_{k})^{k-1}.
$$
By definition of $R_{k}$, we have $(R_{k})^{k} - (R_{k})^{k-1} + 2R_{k} - 1=0$.
It then follows
$$
\frac{F(T-v)}{F(T)} =  \frac{1 + \gamma}{2 + \gamma} < 
\frac{1 + (R_{k})^{k-1}}{2 + (R_{k})^{k-1}} = R_{k},
$$
as desired.
\end{proof}

Now, we may prove Theorem~\ref{th-task2}.

\begin{proof}[Proof of Theorem~\ref{th-task2}]
Let $T$ be an extremal tree. We know by Theorem ~\ref{th-task1} that $T$ is a {\ts} or an odd path. 
Arguing by contradiction, we assume that $T$ is a {\ts} which is not balanced. 
We will show that $T$ admits a good rotation, which contradicts the fact that $T$ is extremal.

The tree $T$ has at least two centers, as otherwise $T$ is trivially balanced.
Let $x$ and $x'$ be distinct centers of $T$ maximizing the difference $\deg(x) - \deg(x')$.
Let also $k:=\deg(x)$ and $\ell:=\deg(x')$ (hence, $k \geq \ell + 2$).

Choose a neighbor $y$ of $x$ that is not on the unique $xx'$-path in $T$, and consider
the rotation $\rho:=(yx,yx')$ (see Figure~\ref{fig-rot-th-task2} for an illustration). Observe that $\rho(T)$ is also a {\ts}, and that $T$ and $\rho(T)$ have the same set of centers. In particular, $\alpha(\rho(T))=\alpha(T)$ by Lemma~\ref{lem-mss}.
Hence, to prove that $\rho$ is a good rotation, it remains to show that $F(\rho(T)) < F(T)$.
This will be done by combining Lemma~\ref{lem-rotating-edge} 
with Lemmas~\ref{lem-Rk} and ~\ref{lem-ToS-leaf-gen}.

\begin{figure}[htb]
\begin{center}
\includegraphics[scale=0.9]{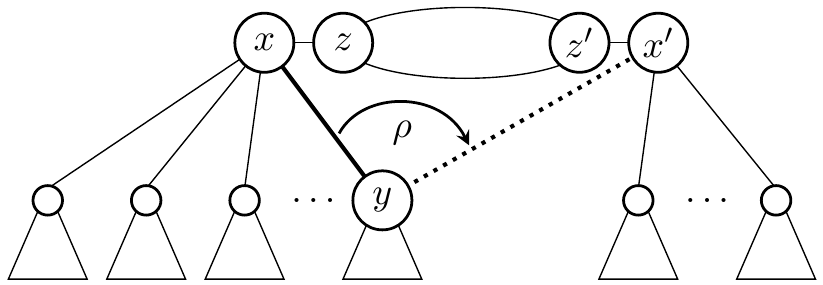}
\caption{\label{fig-rot-th-task2} 
The rotation $\rho$ used in the proof of Theorem~\ref{th-task2}.
}
\end{center}
\end{figure}

By our choice of $x$ and $x'$, every component of $T- \{x , x'\}$ is a {\ts}, every center 
of which has degree at least $\ell$. 
Consider Notation~\ref{nota-rot} with respect to the rotation $\rho=(yx,yx')$. 
(Thus, $z$ and $z'$ are the neighbors of respectively $x$ and $x'$ that lie on the $xx'$-path in $T$, 
and possibly $z=z'$.)
It follows from Lemma~\ref{lem-ToS-leaf-gen} that
$$
\frac{\overline Z_{z}}{Z} < R_{\ell}
$$
and
$$
\frac{\overline X}{X} < (R_{\ell})^{k-2} \leq  (R_{\ell})^{\ell}.
$$
Also, we have that
$$
\frac{\overline Z_{z'}}{Z} \geq \frac{1}{2}
$$
and
$$
\frac{\overline X'}{X'} \geq \frac{1}{2^{\ell-1}},
$$
since $F(T') \leq 2F(T' - v)$ 
holds for every tree $T'$ and vertex $v$ of $T'$ (cf.\ Lemma~\ref{lem:basic}). 
Combining these inequalities with Lemma~\ref{lem-Rk}, we obtain
\begin{align*}
X\overline X' \overline Z_{z'} &> \frac{\overline X}{(R_{\ell})^{\ell}} 
\cdot \frac{X'}{2^{\ell-1}} \cdot  \frac{\overline Z_{z}}{2R_{\ell}} \\
&= \overline X X' \overline Z_{z}\cdot \frac{(R_{\ell})^{-(\ell+1)}}{2^{\ell}} \\
&> \overline X X' \overline Z_{z}\cdot \frac{ \left ( 2^{-\ell/(\ell+1)} \right)^{-(\ell+1)}}{2^{\ell}} \\
&=\overline X X' \overline Z_{z}.
\end{align*}
Therefore, $F(\rho(T)) < F(T)$ by Lemma~\ref{lem-rotating-edge}. This completes the proof.
\end{proof}

\section{Further Results}
\label{sec-task3}

In Section~\ref{sec-task2}, it has been shown that every extremal tree $T$ that is not an odd path must be a balanced {\ts}. An open problem is to understand how these stars are linked together. A few results in this direction are given here, but the general problem is far from being solved. 

\begin{figure}[!htb]
\begin{center}
\includegraphics[scale=0.85]{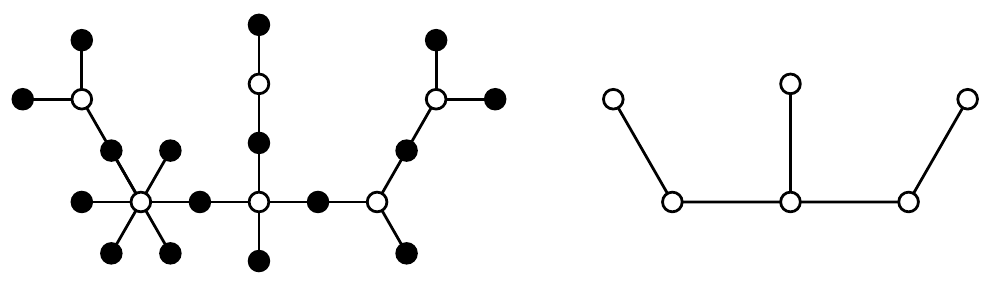}
\caption{\label{fig-ct} 
A {\tos} (left) and its corresponding {\ctree} (right).}
\end{center}
\end{figure}

Let $T$ be a {\ts}. The \emph{{\ctree}} of $T$ is the tree $C_T$ having the centers of $T$ as vertex set, and where two centers are adjacent
if they are at distance 2 in $T$ (see Figure~\ref{fig-ct} for an illustration).

Let $T$ be a balanced {\tos} with $n$ vertices and stability number $\alpha$.
Since $T$ is balanced, every center of $T$ has degree
either $\left\lceil \frac{n-1}{n-\alpha} \right\rceil$ 
or $\left \lceil \frac{n-1}{n-\alpha} \right \rceil - 1$. 
A center is said to be {\DEF heavy} 
in the first case, {\DEF light} in the second.
It can be checked that the number of heavy centers is given by
$$
h(n,\alpha) := 
\left\{
\begin{array}{lll}
(n-1) \mod (n-\alpha)  &   &  \text{if $(n-1) \mod (n-\alpha) \neq 0$} \\
n-\alpha  &   &  \text{otherwise}. \\
\end{array}
\right.
$$
The number of light centers of $T$ is thus
$$
\l(n,\alpha) := n - \alpha - h(n,\alpha).
$$

The next theorem gives a simple characterization of the extremal trees when 
$\l(n,\alpha) \leq 2$.
Let us remark that, if $\alpha = n/2$, then the path $P_{n}$
is the only tree with $n$ vertices
and stability number $\alpha$ that is extremal, as follows from   
the result of Prodinger and Tichy~\cite{Prodinger82} mentioned in the introduction. 
Hence, we assume $\alpha > n/2$ in what follows.

\begin{thm}
\label{th-ct-path}
Let $T$ be a tree with $n$ vertices and stability number $\alpha > n/2$, and
assume $\l(n,\alpha) \leq 2$. Then, $T$ is extremal if and only if 
$T$ satisfies the following three conditions:
\begin{itemize}
\item $T$ is a balanced {\ts};
\item the {\ctree} of $T$ is isomorphic to a path $P$, and
\item each light center of $T$ (if any) is an endpoint of $P$.
\end{itemize}
\end{thm}

Our proof of Theorem~\ref{th-ct-path} is based on the following lemma.

\begin{lem} \label{lem-task3-a}
Let $T$ be a balanced {\ts}, and assume that a leaf $w$ of the {\ctree} $C_{T}$ is heavy.
If $T$ is extremal, then $C_{T}$ is a path, and all internal vertices of $C_{T}$ are heavy.
\end{lem}
\begin{proof}
Arguing by contradiction, suppose that $T$ is extremal and contains a center $v$
such that either $v$ is heavy with degree at least $3$ in $C_{T}$, or light with degree at least $2$ in $C_{T}$.
We may assume that $v$ has been chosen so that every inner vertex of the $vw$-path in $C_{T}$ 
is heavy and has degree $2$ in $C_{T}$.

Let $P$ be the unique $vw$-path in the tree $T$. 
Let $y$ be a neighbor of $v$ in the {\ctree} $C_{T}$ such that $y \notin V(P)$. Let $x$ be the vertex of $T$
that is adjacent to both $v$ and $y$. 
Let $x'$ be a leaf of $T$ that is adjacent to $w$. (Thus, $x' \notin V(P)$.)
Let $\rho$ be the rotation $\rho:=(yx,yx')$ (see Figure~\ref{fig-rot-lem-task3-a} for an illustration).
The tree $\rho(T)$ is a {\ts} and has the same number of centers as $T$, thus
$\alpha(\rho(T))=\alpha(T)$ by Lemma~\ref{lem-mss}. Hence, to reach a contradiction,
it is enough to show that $F(\rho(T)) < F(T)$.

\begin{figure}[htb]
\begin{center}
\includegraphics[scale=0.9]{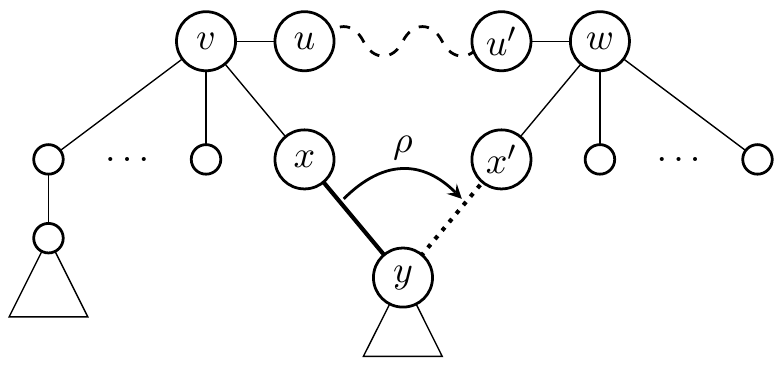}
\caption{\label{fig-rot-lem-task3-a} 
The rotation $\rho$ used in the proof of Lemma~\ref{lem-task3-a}.
}
\end{center}
\end{figure}

Consider Notation~\ref{nota-rot} with respect to the rotation $\rho$.
Thus $z = v$, $z' = w$, and $X = \overline X = X' = \overline X' = 1$.
By Lemma~\ref{lem-rotating-edge}, $F(\rho(T)) < F(T)$ if and only if $\overline Z_{z} < \overline Z_{z'}$. 

Let $u$ and $u'$ be the neighbors of respectively $z$ and $z'$ in the path $P$.
Let $\widetilde T$ be the component of $T - \{z, z'\}$ that contains both $u$ and $u'$.
Since every center of $T$ included in $V(P) - \{z, z'\}$ is heavy and has degree $2$ in $C_{T}$, the trees
$\widetilde T - u$ and $\widetilde T - u'$ are isomorphic. 

Let $v_{1}, \dots, v_{k}$ be the neighbors of $z$ ($=v$) in $T$ that are distinct from $u$ and $x$.
Let $T^{z}_{i}$ ($1 \leq i \leq k$) be the component of $T - \{z\}$ that includes $v_{i}$.
Let $\l:= \deg_{T}(z') - 2$. Thus, $\l = k$ or $k+1$, depending on whether $z$ is heavy or light.
Since the $\l$ neighbors of $z'$ that are distinct from $u'$ and $x'$ are all leaves of $T$, it follows that
\begin{align*}
\overline Z_{z} &=  \left( \prod_{i=1}^k F(T^{z}_{i}) \right) \cdot \left(2^{\l}F(\widetilde T) + F(\widetilde T - u')  \right), \\
\overline Z_{z'} &=  2^{\l} \cdot \left(F(\widetilde T)\prod_{i=1}^k F(T^{z}_{i}) +  F(\widetilde T - u)\prod_{i=1}^k F(T^{z}_{i} - v_{i})  \right).
\end{align*}
Using that $F(\widetilde T - u) = F(\widetilde T - u')$, we obtain
$$
\overline Z_{z'} - \overline Z_{z} = F(\widetilde T - u) \cdot \left(2^{\l}\prod_{i=1}^k F(T^{z}_{i} - v_{i}) - \prod_{i=1}^k F(T^{z}_{i})\right).
$$
We have $F(T^{z}_{i}) \leq 2F(T^{z}_{i} - v_{i})$ for every $i \in \{1, \dots, k\}$, with strict inequality if $|T^{z}_{i}| > 1$.
If $z$ is light, then $\l = k + 1$, and
$$
2^{\l}\prod_{i=1}^k F(T^{z}_{i} - v_{i}) - \prod_{i=1}^k F(T^{z}_{i}) \geq
(2^{\l} - 2^{k})\prod_{i=1}^k F(T^{z}_{i} - v_{i}) > 0.
$$
Similarly, if $z$ is heavy, then $\l = k$. Moreover, $|T^{z}_{j}| > 1$ holds for some $j \in \{1, \dots, k\}$,
because $z$ has degree at least $3$ in $C_{T}$. This implies $F(T^{z}_{j}) < 2F(T^{z}_{j} - v_{j})$, and
$$
2^{\l}\prod_{i=1}^k F(T^{z}_{i} - v_{i}) - \prod_{i=1}^k F(T^{z}_{i}) >
(2^{\l} - 2^{k})\prod_{i=1}^k F(T^{z}_{i} - v_{i}) = 0.
$$
Thus, $\overline Z_{z} < \overline Z_{z'}$ holds in both cases, as claimed. This concludes the proof.
\end{proof}

\begin{proof}[Proof of Theorem~\ref{th-ct-path}]
First, we observe that there is a unique tree (up to isomorphism) satisfying the three conditions 
given in the statement of the theorem. Hence, it is enough to show that each of these three conditions
is necessary for $T$ to be extremal.

Thus, suppose that $T$ is extremal.
Then $T$ is a balanced {\ts} by Theorem~\ref{th-task2} 
(note that $T$ cannot be an odd path since $\alpha > n/2$).
If the {\ctree} of $T$ is not a path, then it has
at least three leaves, and one of them is heavy.
But then Lemma~\ref{lem-task3-a} implies that $T$ is not extremal, a contradiction. 
Hence, the {\ctree} is isomorphic to a path $P$. Furthermore, all internal vertices
of $P$ are heavy centers of $T$, by the same lemma. The theorem follows.
\end{proof}

\begin{figure}
\begin{center}
\includegraphics[scale=0.9]{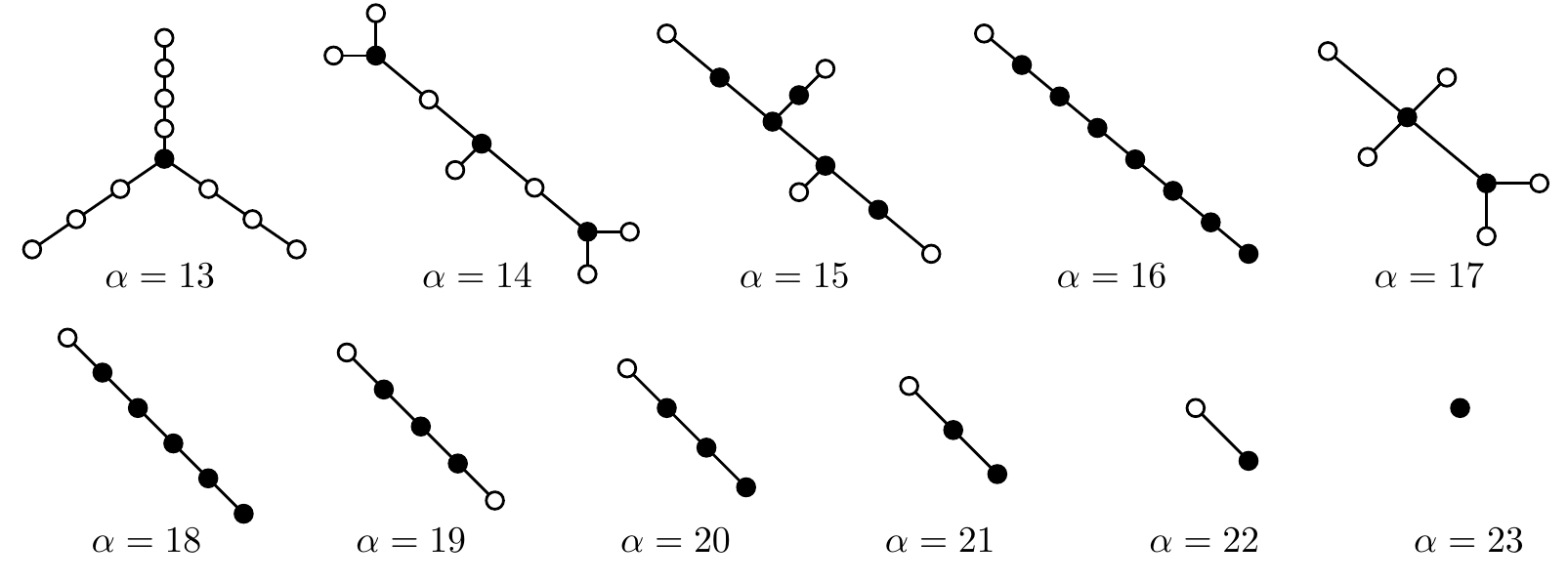}
\caption{\label{fig-center-trees-n-24}
The {\ctree s} of all {\tss} on $n=24$ vertices which are extremal.
Light centers are drawn in white and heavy centers in black.
}
\end{center}
\end{figure}

When $\l(n,\alpha) \geq 3$ and $\alpha > n/2$, every extremal tree $T$ is 
a balanced {\ts} by Theorem~\ref{th-task2}. However, the {\ctree} of $T$
is no longer necessarily a path.
This can be seen on Figure~\ref{fig-center-trees-n-24}, which provides a list
of the {\ctree s} of all extremal trees for $n=24$ and $\alpha > 12$ (these have been computed with 
the system GraPHedron~\cite{Melot08}).

\section*{Acknowledgments}
The authors thank Guy Louchard for useful discussions, and 
the two anonymous referees for their helpful comments on a previous version of the manuscript. 
This work was supported by the {\em Communaut\'e Fran\c caise de Belgique (projet ARC)}.

\bibliography{Fibo}
\bibliographystyle{plain}

\end{document}